\documentclass[12pt,reqno]{amsart}
\usepackage{amssymb}

\let\frak\mathfrak

\usepackage{eucal}

\def\>{\relax\ifmmode\mskip.666667\thinmuskip\relax\else\kern.111111em\fi}
\def\<{\relax\ifmmode\mskip-.333333\thinmuskip\relax\else\kern-.0555556em\fi}
\def\vsk#1>{\vskip#1\baselineskip}
\def\vv#1>{\vadjust{\vsk#1>}\ignorespaces}
\def\vvn#1>{\vadjust{\nobreak\vsk#1>\nobreak}\ignorespaces}

\def\plait#1{\par\hangindent2\parindent\indent\kern\parindent
\llap{#1\enspace}\ignorespaces}

\let\Smallskip\smallskip
\def\smallskip{\par\Smallskip}
\let\Medskip\medskip
\def\medskip{\par\Medskip}
\let\Bigskip\bigskip
\def\bigskip{\par\Bigskip}

\let\Maketitle\maketitle
\def\maketitle{\Maketitle\thispagestyle{empty}\let\maketitle\empty}

\newtheorem{theorem}{Theorem}[section]

\newtheorem{cor}[theorem]{Corollary}
\newtheorem{lem}[theorem]{Lemma}

\newtheorem{thm}[theorem]{Theorem}
\theoremstyle{definition}

\theoremstyle{remark}

\theoremstyle{remark}

\numberwithin{equation}{section}

\let\leq\leqslant
\let\ge\geqslant
\let\geq\geqslant

\newcommand{\nc}{\newcommand}
\nc{\on}{\operatorname}
\nc{\ch}{\mbox{ch}}
\nc{\Z}{{\mathbb Z}}
\nc{\C}{{\mathbb C}}
\nc{\R}{{\mathbb R}}
\nc{\pone}{{\mathbb C}{\mathbb P}^1}
\nc{\pa}{\partial}
\nc{\F}{{\mathcal F}}
\nc{\arr}{\rightarrow}
\nc{\larr}{\longrightarrow}
\nc{\al}{\alpha}
\nc{\ri}{\rangle}
\nc{\lef}{\langle}
\nc{\W}{{\mathcal W}}
\nc{\la}{\lambda}
\nc{\ep}{\epsilon}
\nc{\eps}{\varepsilon}

\nc{\Om}{\Omega}
\nc{\su}{\widehat{{\mathfrak sl}}_2}
\nc{\sw}{{\mathfrak s}{\mathfrak l}}

\nc{\g}{{\mathfrak g}}
\nc{\h}{{\mathfrak h}}
\nc{\n}{{\mathfrak n}}
\nc{\N}{\widehat{\n}}
\nc{\G}{\widehat{\g}}
\nc{\De}{\Delta_+}
\nc{\gt}{\widetilde{\g}}
\nc{\Ga}{\Gamma}
\nc{\one}{{\mathbf 1}}
\nc{\z}{{\mathfrak Z}}
\nc{\zz}{{\mathcal Z}}
\nc{\Hh}{{\mathcal H}_\beta}
\nc{\qp}{q^{\frac{k}{2}}}
\nc{\qm}{q^{-\frac{k}{2}}}
\nc{\La}{\Lambda}
\nc{\wt}{\widetilde}
\nc{\qn}{\frac{[m]_q^2}{[2m]_q}}
\nc{\cri}{_{\on{cr}}}
\nc{\kk}{h^\vee}
\nc{\sun}{\widehat{\sw}_N}
\nc{\hh}{\widehat{\mathfrak h}}
\nc{\HH}{{\mathcal H}_{q,t}}
\nc{\ca}{\wt{{\mathcal A}}_{h,k}(\sw_2)}
\nc{\gl}{\widehat{{\mathfrak g}{\mathfrak l}}_2}
\nc{\el}{\ell}
\nc{\s}{{\mathbf s}}
\nc{\bi}{\bibitem}
\nc{\om}{\omega}
\nc{\WW}{\W_\beta}
\nc{\scr}{{\mathbf S}}
\nc{\ab}{{\mathbf a}}
\nc{\rr}{r}
\nc{\ol}{\overline}
\nc{\con}{qt^{-1} + q^{-1}t}
\nc{\den}{q^{\el-1} t^{-\el+1}+ q^{-\el+1} t^{\el-1}}
\nc{\ds}{\displaystyle}
\nc{\B}{\mc B}
\nc{\A}{{\mathbb A}}
\nc{\GG}{{\mathcal G}}
\nc{\UU}{{\mathcal U}}
\nc{\MM}{{\mathcal M}}
\nc{\CC}{{\mathcal C}}
\nc{\GL}{{}^L G}
\nc{\dzz}{\frac{dz}{z}}
\nc{\Res}{\on{Res}}
\nc{\rep}{{\mathcal R}ep \;}
\nc{\uqg}{U_q \G}
\nc{\uqgg}{U_q \g}
\nc{\Fq}{{\mathbb F}_q}

\nc{\stimes}{\ltimes}
\nc{\K}{\hat{\mathcal K}}
\nc{\Ql}{\ol{\mathbb Q}_\ell}
\renewcommand{\O}{\hat{\mathcal O}}
\nc{\ga}{\gamma}
\nc{\PL}{{}^L P}
\nc{\E}{\mc E}
\nc{\mc}{\mathcal}
\nc{\mbf}{\mathbf}
\nc{\bb}{{\mathfrak b}}
\nc{\OO}{{\mc O}}
\nc{\Po}{{\mc P}}
\nc{\V}{{\mc V}}
\nc{\yy}{{\mc Y}}
\nc{\M}{\mathcal M}
\nc{\Coh}{{{\mathcal C}oh}}
\nc{\Cohn}{\Coh_n}
\nc{\f}{{\mathcal F}}
\nc{\si}{_E}
\nc{\Gaf}{{\mathbb G}_{a,\Fq}}
\nc{\KK}{{\mathfrak k}}

\nc{\PO}{{\mathbb P^1}}
\nc{\PR}{{\mathbb P^r}}

\nc{\Wr}{{ {\rm Wr}}}

\newcommand{\bean}{\begin{eqnarray}}
\newcommand{\eean}{\end{eqnarray}}
\newcommand{\be}{\begin{displaymath}}
\newcommand{\ee}{\end{displaymath}}
\newcommand{\bea}{\begin{eqnarray*}}
\newcommand{\eea}{\end{eqnarray*}}
\newcommand{\beq}{\begin{equation}}
\newcommand{\eeq}{\end{equation}}
\newcommand{\bs}{\boldsymbol}
\newcommand{\Ref}[1]{{$($\ref{#1}$)$}}

\newcommand{\sing}{{\rm Sing}\,}

\newcommand{\glN}{{\frak{gl}_{N+1}}}
\newcommand{\D}{{\mathcal D}}
\newcommand{\End}{{\rm End}}
\newcommand{\T}{{\bs p}}
\newcommand{\nash}{L_{\bs \La}[\bs\la^{(\infty)}]}
\newcommand{\snash}{\sing L_{\bs \La}[\bs\la^{(\infty)}]}

\newcommand{\AT}{ A_{\T,\Phi}}

\def\gln{\frak{gl}_{N+1}}
\def\glnt{\gln[t]}

\def\Uglnt{U(\glnt)}

\def\rdet{\on{rdet}}
\def\Gr{{\rm Gr}}
\def\O{{\mc O}}

\let\der\partial

\def\lab{{{{\bs\la^{(\vee)}}}}}
\def\Lbl{{\bs\La,\bs\la^{(\infty)}\!,\bs z}}

\def\lc{,\dots,}
\def\dl{\delta}
\def\si{\sigma}

\let\kk K

\nc {\LL}{L_{\bs\La}}
\nc{\AB}{{A_{\T,B}}}
\nc{\BL}{{\B_{\bs\La,\bs\la^{(\infty)}\!,\bs z}}}

\begin{document}

\title[Bethe algebra and algebra of functions on critical points]
{Bethe algebra of the $\glN$ Gaudin model and\\[2pt] algebra of functions
on the critical set\\[2pt] of the master function}

\author[E.\,Mukhin, V\<.\,Tarasov, and A.\>Varchenko]
{E.\,Mukhin ${}^{*,1}$, V\<.\,Tarasov ${}^{*,\star,2}$,
A.\>Varchenko {${}^{**,3}$} }

\thanks{${}^1$\ Supported in part by NSF grant DMS-0900984}
\thanks{${}^2$\ Supported in part by NSF grant DMS-0901616}
\thanks{${}^3$\ Supported in part by NSF grant DMS-0555327}

\maketitle

\centerline{\it ${}^{*}$Department of Mathematical Sciences,
Indiana University -- Purdue University,}
\centerline{\it Indianapolis, 402 North Blackford St, Indianapolis,
IN 46202-3216, USA}
\smallskip
\centerline{\it $^\star$St.\,Petersburg Branch of Steklov Mathematical
Institute}
\centerline{\it Fontanka 27, St.\,Petersburg, 191023, Russia}
\smallskip
\centerline{\it ${}^{**}$Department of Mathematics, University of
North Carolina at Chapel Hill,} \centerline{\it Chapel Hill, NC
27599-3250, USA}
\bigskip

\begin{abstract}
Consider a tensor product of finite-dimensional irreducible $\glN$-modules and
its decomposition into irreducible modules. The $\glN$ Gaudin model assigns to
each multiplicity space of that decomposition a commutative (Bethe) algebra of
linear operators acting on the multiplicity space. The Bethe ansatz method is
a method to find eigenvectors and eigenvalues of the Bethe algebra. One starts
with a critical point of a suitable (master) function and constructs
an eigenvector of the Bethe algebra.

In this paper we consider the algebra of functions on the critical set of
the associated master function and show that the action of this algebra on
itself is isomorphic to the action of the Bethe algebra on a suitable subspace
of the multiplicity space.

As a byproduct we prove that the Bethe vectors corresponding to different
critical points of the master function are linearly independent and, in
particular, nonzero.

\end{abstract}

\section{Introduction}

Let $L_{\bs\la}$ be the irreducible finite-dimensional $\glN$-module of highest
weight $\bs\la$. Let $L_{\bs\La}=\otimes_{s=1}^n L_{\bs\la^{(s)}}$ be a tensor
product of such modules, and
$L_{\bs \La}=\oplus_{\bs\la^{(\infty)}}L_{\bs\la^{(\infty)}}\otimes
W_{\bs\la^{(\infty)}}$ the decomposition into
irreducible representations. The multiplicity space
$W_{\bs\la^{(\infty)}}$ of $L_{\bs\la^{(\infty)}}$
can be identified with $\snash \subset L_{\bs\La}$, the subspace of
singular vectors of weight $\bs\la^{(\infty)}$.
To each multiplicity space $\snash$ and distinct complex numbers
$z_1,\dots,z_n$, the $\glN$ Gaudin
model assigns a
commutative subalgebra of $\End\,(\snash)$ called the Bethe algebra
and denoted by $\BL$.

The $\BL$-module $\snash$ has an interesting geometric realization.
In \cite{MTV3}
we constructed an isomorphism of the
$\BL$-module $\snash$
and the regular representation of the algebra of functions on the
scheme-theoretical intersection of suitable Schubert cycles.
This isomorphism
can be viewed as the geometric
Langlands correspondence in the $\glN$ Gaudin model.
In [MTV4] we argued that this geometric Langlands correspondence extends
to the third, equally important, player --- the algebra of functions on
the critical set of the corresponding master function.
In this paper we prove another result supporting that principle.

The master and weight functions are useful objects associated with each multiplicity space
$\snash$,
see \cite{SV}. They are functions of some auxiliary variables
$\bs t=(t^{(i)}_j)$. The master function $\Phi(\bs t)$ is a scalar function
and the weight function $\omega(\bs t)$ is an $L_{\bs\La}$-valued function.
They are used in the Bethe ansatz method to construct eigenvectors of the Bethe algebra
$\BL$. Namely, if $\T$ is a critical point of the master function then
the vector $\omega(\T)$ lies in $\snash$ and is an eigenvector
of the Bethe algebra, see \cite{MTV2}.

In this paper we consider the algebra $A_\Phi$ of functions on the
critical set
of the master function. With the help of the weight
function, we construct a linear embedding $\al : A_\Phi \to \snash$
and show that $\al(A_\Phi)$ is a $\BL$-submodule of $\snash$.
We denote the image of $\B$ in $\End(\al(A_\Phi))$ by $A_B$. We construct an algebra
isomorphism $\beta : A_\Phi \to A_B$ and show that the $A_B$-module
$\al(A_\Phi)$ is isomorphic to the regular representation of $A_\Phi$.
That statement is our main result, see Theorems \ref{thm main} and \ref{thm many
point}. As a byproduct
we show that for
any critical point $\T$ of the master function, the vector
$\omega(\bs \T)$ is nonzero. For a nondegenerate critical point
that fact was proved in \cite{MV2} and \cite{V}.

The paper is organized as follows. In Section \ref{Algebra APhi} we define
the master function $\Phi$ and the algebra $A_\Phi$
of functions on the critical set
$C_\Phi$ of the master function. The algebra $A_\Phi$ is the direct sum of local algebras
$\AT$ corresponding to points $\T\in C_\Phi$. In Theorem
\ref{thm generate}
we describe useful generators of the algebra $\AT$.
We prove Theorem \ref{thm generate} in Section \ref{Algebra AGr}.
In Section \ref{sec algebra A_B} the Bethe algebra is introduced.
We define the weight function in Section \ref{sec weighT} and formulate our first main result
Theorem \ref{thm main}.
We prove Theorem \ref{thm main} in Section \ref{pRoofs}. Our second main result,
Theorem \ref{thm many point}
is formulated and proved in Section
\ref{Concluding remarks}.

\medskip
The authors thank P.\,Belkale for help in proving Theorem~\ref{thm generate}.

\section{Algebra $A_\Phi$}
\label{Algebra APhi}

\subsection{Lie algebra $\gln$}

Let $e_{ij}$, $i,j=1\lc N+1$, be the standard generators of the Lie algebra
$\gln$ satisfying the relations
$[e_{ij},e_{sk}]=\dl_{js}e_{ik}-\dl_{ik}e_{sj}$. Let $\h\subset\gln$ be
the Cartan subalgebra generated by $e_{ii}, \,i=1\lc N+1$.
Let $\h^*$ be the dual space. Let
$\epsilon_i, \,i=1\lc N+1,$ be the basis of $\h^*$ dual to the basis
$e_{ii}, \,i=1\lc N+1$, of $\h$.
Let $\al_1,\dots,\al_{N} \in \h^*$ be simple roots, $\al_i
=\epsilon_i-\epsilon_{i+1}$. Let
$(\,,\,)$ be the standard scalar product on $\h^*$ such that
the basis $\epsilon_i, \,i=1\lc N+1,$
is orthonormal.

A sequence of integers $\bs\la=(\la_1,\dots,\la_{N+1})$ such that
$\la_1\ge\la_2\ge\dots\ge\la_{N+1}\ge0$ is called a {\it partition with at most
${N+1}$ parts\/}. Denote $|\bs\la|=\sum_{i=1}^{N+1}\la_i$.
We identify partitions $\bs\la$
with vectors $\la_1\epsilon_1+
\dots + \la_{N+1}\epsilon_{N+1}$ of $\h^*$.

\subsection{Master function}
\label{subsubsection on master function}

Let $\bs \La = (\bs\la^{(1)},\dots, \bs\la^{(n)})$ be a collection of partitions,
where $\bs\la^{(i)}=(\la^{(i)}_1,\dots,\la^{(i)}_{N+1})$ and
$\la^{(i)}_{N+1}=0$.
Let $\bs l = (l_1,\dots,l_N)$ be nonnegative integers such that
$$
\bs\lambda^{(\infty)}\ =\ \sum_{i=1}^n \bs\la^{(i)} - \sum_{j=1}^{N} l_j\al_j
$$
is a partition.
Denote $l=l_1+\dots +l_N$,
\be
\bs t\ =\ (t^{(1)}_1,\dots,t^{(1)}_{l_1},t^{(2)}_{1},\dots,t^{(2)}_{l_2},\dots,
t^{(N)}_{1},\dots,t^{(N)}_{l_{N}})\ .
\ee
Fix a collection of distinct complex numbers
$\bs z = (z_1,\dots,z_n)$. Let $\Phi (\bs t)$ be
the master function associated with this data,
\be
\Phi (\bs t) =
\prod_{i=1}^{N}\prod_{1\leq j<j'\leq l_i} \!\!\!(t_j^{(i)}-t_{j'}^{(i)})^{2}
\prod_{i=1}^{N-1}
\prod_{j=1}^{l_i}\prod_{j'=1}^{l_{i+1}}
(t_j^{(i)}-t_{j'}^{(i+1)})^{-1}
\prod_{i=1}^{N}\prod_{j=1}^{l_{i}}
\prod_{s=1}^{n}
(t_{j}^{(i)}-z_s)^{-(\bs\la^{(s)},\al_i)} .
\ee
Denote
\beq
\label{def of U}
U = \{ \T \in \C^l\ | \ \Phi
\ \text{is well-defined at}\ \T\ {\rm and}\ \Phi(\T)\neq 0\} .
\eeq
The set
$U$ is the complement in $\C^l$ to a union of hyperplanes.
The master function is a rational function regular on $U$.

Denote by
$\C(\bs t)_U$
the algebra of rational functions on $\C^l$ regular on $U$.
The partial derivatives
\be
\Psi_{ij}\ =\ \partial(\log \Phi)/\partial t^{(i)}_j ,
\qquad
i=1,\dots, N, \quad j=1,\dots,l_i\ ,
\ee
are elements of $\C(\bs t)_U$.
Denote by $I_\Phi\subset \C(\bs t)_U$ the ideal generated by
$\Psi_{ij},$ $i=1,\dots, N,$ $ j=1,\dots,l_i$, and set
\beq
\label{APhi}
A_\Phi\ =\ \C(\bs t)_U/ I_\Phi .
\eeq
Denote by $C_\Phi$ the zero set of the ideal. The zero set is
finite, \cite{MV1}. The algebra $A_\Phi$ is finite-dimensional and
is the direct sum of local algebras,
\be
A_\Phi\ =\ \oplus_{\T\in C_\Phi} A_{\T, \Phi}\
\ee
corresponding to points $\T\in C_\Phi$. For $ {\T}\in C_\Phi$,
the local algebra $A_{\T, \Phi}$ may be defined as the
quotient of the algebra of germs at $\T$ of holomorphic functions
modulo the ideal $I_{\T,\Phi}$ generated by all the functions
$\Psi_{ij}$.
The algebra $\AT$ contains the
maximal ideal $\frak m_\T$ generated by the germs
of functions equal to zero at $\T$.

\subsection{Generators of the local algebra of a critical point}
\label{sec generators}
Let $u$ be a variable.
Define an $N$-tuple of polynomials $T_1,\dots, T_N\in \C[u]$,
\be
T_i(u)\ =\ \prod_{s=1}^n\,(u-z_s)^{(\bs\la^{(s)},\al_i)}\ ,
\ee
an $N$-tuple of polynomials
$y_1,\dots,y_N\in \C[u,\bs t]$,
\be
y_i(u, \bs t)\ =\ \prod_{j=1}^{l_i}\, (u- t^{(i)}_j)\ ,
\ee
and the differential operator
\begin{align*}
\D_{\Phi}\, &{}=\,
(\der_u - \log' ( \frac { T_1\dots T_N } { y_{1} } ) )
\\[2pt]
&\<\>{}\times\,( \der_u - \log' ( \frac {y_{1}T_2\dots T_{N} } {y_{2} } ) )
\dots ( \der_u - \log' ( \frac {y_{N-1} T_N}{ y_N } ) )
( \der_u - \log' ( y_N ) ) ,
\end{align*}
where $\der_u = d/du$ and $\log' f$ denotes $ (df/du)/f$. We have
\beq
\label{formula h(t,u)}
\D_{\Phi}
\ =\
\der_u^{N+1}+\sum_{i=1}^{N+1}\,G_i\,\der_u^{N+1-i},
\qquad
G_i = \sum_{j=i}^\infty G_{ij} u^{-j} ,
\eeq
where $G_{ij} \in \C[\bs t]$.

For $\T\in C_\Phi$, $f\in \C(\bs t)_U$ denote by $\bar f$ the
image of $f$ in $\AT$.
Denote
\be
\bar\D_{\Phi}
\ =\
\der_u^{N+1}+\sum_{i=1}^{N+1}\,\bar G_i\,\der_u^{N+1-i},
\ee
where
$\bar G_i = \sum_{j=i}^\infty \bar G_{ij} u^{-j}$.

\begin{thm}
\label{thm generate}
For any $\T\in C_\Phi$,
the elements
$\bar G_{ij},$ $i=1,\dots,N$, $j\geq i$, generate $\AT$.
\end{thm}

Theorem \ref{thm generate} is proved in Section \ref{PRoofs}.

\subsection{Polynomials $h_i$}
\label{sec 1.3}
Let $A$ be a commutative algebra.
For $g_1,\dots,g_{i} \in A[u]$, denote by
$\Wr(g_1(u),\dots,g_{i}(u))$ the { Wronskian},
\be
\Wr(g_1(u),\dots,g_{i}(u))\,=\,
\det\left(\begin{matrix} g_1(u) & g_1'(u) &\dots & g_1^{(i-1)}(u) \\
g_2(u) & g_2'(u) &\dots & g_2^{(i-1)}(u) \\ \dots & \dots &\dots & \dots \\
g_{i}(u) & g_{i}'(u) &\dots & g_{i}^{(i-1)}(u)
\end{matrix}\right),
\ee
where $g^{(j)}(u)$ denotes the $j$-th derivative of $g(u)$ with respect to $u$.

Introduce a set
\beq
\label{set P}
P\,=\,\{d_1,d_2,\dots,d_{N+1}\}\,,\qquad d_i=\la^{(\infty)}_i+{N+1}-i\,.
\eeq

\begin{thm}
\label{thm Bmv}
There exist unique polynomials $h_1,\dots,h_{N+1} \in \AT[u]$ of the form
\beq
\label{h basis}
h_i=u^{d_i}+\sum_{j=1,\ d_i-j\not\in P}^{d_i}h_{ij}u^{d_i-j}
\eeq
such that $h_{N+1}=y_{N}$ and
\beq
\label{w h}
\Wr(h_{N+1},h_{N}\dots,h_{N+1-j})\ = \
y_{N-j}T^{j}_NT^{j-1}_{N-1}\dots T^{1}_{N-j+1}
\!\!\!\!\!
\prod_{N+1-j\leq i<i'\leq N+1}
\!\!\!\!\!
(d_{i}-d_{i'})\
\eeq
for $j=1,\dots,N$, where $y_0=1$. Moreover, each of the polynomials
$h_1,\dots,h_{N+1} $ is a solution of the differential equation
$\bar \D_\Phi h(u) = 0$.
\qed
\end{thm}

\begin{proof}
The existence of unique polynomials $h_i$ satisfying \Ref{w h} is proved in
~\cite{BMV} generalizing the corresponding result in Section~5 of~\cite{MV1}.
The fact that the polynomials $h_i$ satisfy the differential equation
$\bar \D_\Phi h(u) = 0$ is proved like in Section~5 of~\cite{MV1}.
\end{proof}

\begin{lem}
\label{lem Generate}
The subalgebra of $\AT$ generated by elements
$\bar G_{ij},$ $i=1,\dots,N$, $j\geq i$, contains all the coefficients
$h_{ij} ,\ i=1,\dots,{N+1},\ j=1,\dots,d_i,\ d_i-j\not\in P$.
\end{lem}

The proof of the lemma is the same as the proof of Lemma 3.4 in \cite{MTV3}.

\section{Algebra $A_{\Gr}$}
\label{Algebra AGr}

\subsection{Algebra
$\O_{\bs\la^{(\infty)}}$}
\label{Ominfty}
Let $\bs \La = (\bs\la^{(1)},\dots, \bs\la^{(n)})$, $\bs\la^{(\infty)}$,
$\bs z=(z_1,\dots,z_n)$ be partitions and numbers
as in Section \ref{subsubsection on master function}.
Let $d$ be a natural number such that $d-N-1\geq \la^{(\infty)}_1$
and $d-N-1\geq \la^{(i)}_1$ for $i=1,\dots,n$.

Let $\C_d[u]$ be the space of
polynomials in $u$ of degree less than $d$.
Let $\Gr({N+1},d)$ be the Grassmannian of all ${N+1}$-dimensional subspaces of
$\C_d[u]$.

For a complete flag
$\F=\{0\subset F_1\subset F_2\subset\dots\subset F_d=\C_d[u]\}$ and
a partition $\bs\la=(\la_1,\dots,\la_{N+1})$ with $\la_1\leq d-N-1$, define
a { Schubert cell\/} \,$\Om_{{\bs\la}}(\F)\subset \Gr({N+1},d)$,
\be
\Om_{\bs\la}(\F)=\{\bs q\in\Gr({N+1},d)\ |\ \dim (\bs q\cap F_{d-j-\la_j})={N+1}-j\,,
\ \dim (\bs q\cap F_{d-j-\la_j-1})={N}-j\}\,.
\ee
We have \;$\on{codim}\,\Om_{\bs\la}(\F)=|\bs\la|$.

\medskip

Let $P\,=\,\{d_1,d_2,\dots,d_{N+1}\}$ be defined in \Ref{set P}.
Introduce a new partition
\beq
\label{dual weight}
\lab\,=\,(d-N-1-\la^{(\infty)}_N,d-N-1-\la^{(\infty)}_{N-1},\dots, d-N-1-\la^{(\infty)}_1)\,.
\eeq
Denote
\be
\,\F(\infty)=
\{0\subset \C_1[u]\subset\C_2[u]\subset\dots\subset\C_d[u]\}\,.
\ee
Consider the Schubert cell $\Om_\lab(\F(\infty))$.
We have $\dim\Om_\lab(\F(\infty))=|\bs\la^{(\infty)}|$.

The Schubert cell $\Om_\lab(\F(\infty))$ consists
of ${N+1}$-dimensional subspaces $\bs q\subset\C_d[u]$ with a basis
$\{f_1,\dots,f_{N+1}\}$ of the form
\beq
\label{basis}
f_i=u^{d_i}+\sum_{j=1,\ d_i-j\not\in P}^{d_i}f_{ij}u^{d_i-j}.
\eeq
Such a basis is unique.

Denote by $\O_{\bs\la^{(\infty)}}$ the algebra of regular functions on
$\Om_\lab(\F(\infty))$. The cell $\Om_\lab(\F(\infty))$ is an affine space
with coordinate functions $f_{ij}$. The algebra
$\O_{\bs\la^{(\infty)}}$
is the polynomial algebra in variables $f_{ij}$,
\beq
\label{Ola}
\O_{\bs\la^{(\infty)}}=
\C[f_{ij},\ i=1,\dots,{N+1},\ j=1,\dots,d_i,\ d_i-j\not\in P].
\eeq

\subsection{Intersection of Schubert cells}
\label{Intersection of Scuibert cells}

For $z\in\C$, consider the complete flag
\be
\F(z)\,=\,\bigl\{0\subset (u-z)^{d-1}\C_1[u]\subset(u-z)^{d-2}\C_2[u]
\subset\dots\subset\C_d[u]\bigr\}\,.
\ee
Denote by $\Om_\Lbl$ the set-theoretic intersection and by
$A_\Gr$
the scheme-theoretic intersection of the $n+1$ Schubert cells
$\Om_\lab(\F(\infty))$,
$\Om_{\bs\la^{(s)}}(\F(z_s))\,,s=1,\dots,n$, see Section 4 of
\cite{MTV3}. The set-theoretic intersection
$\Om_\Lbl$
is a finite set
and the scheme-theoretic intersection $A_\Gr$ is a finite-dimensional algebra
(``of functions on $\Om_\Lbl$'').
The algebra
of functions on $\Om_\Lbl$ is the direct sum of local algebras,
\be
A_\Gr\ =\ \oplus_{\bs q \in \Om_\Lbl} A_{\bs q, \Gr}\ ,
\ee
corresponding to points $\bs q\in\Om_\Lbl$.
The algebra $A_\Gr$ is the quotient of the algebra
$\O_{\bs\la^{(\infty)}}$
of functions on $\Om_\lab(\infty)$ by a suitable ideal. For
${\bs q\in\Om_\Lbl}$ and $f\in \O_{\bs\la^{(\infty)}}$ denote by
$\bar f$ the image of $f$ in $A_{\bs q,\Gr}$.

\begin{lem}
\label{lem proJ}
For any $\bs q\in\Om_\Lbl$, the elements
$\bar f_{ij},\ i=1,\dots,{N+1},\ j=1,\dots,d_i,\\ d_i-j\not\in P$,
generate $A_{\bs q,\Gr}$.
\qed
\end{lem}

\subsection{Isomorphism of algebras}
\label{Isomorphism of algebras}

\begin{thm}
[\cite{MV1}]
\label{thm mV}
Let $\T \in C_\Phi$.
Let
$h_1,\dots,h_{N+1} \in A_{\T,\Phi}[u]$ be polynomials defined in Theorem \ref{thm Bmv}.
Denote by $\tilde h_1,\dots,\tilde h_{N+1}$ the projection of the polynomials to
$A_{\T,\Phi}/\frak m_\T [u]=\C[u]$. Then
$\langle\tilde h_{1},\dots,\tilde h_{N+1}\rangle \in \Om_\Lbl$.
\qed
\end{thm}

Denote $\bs q=\langle\tilde h_{1},\dots,\tilde h_{N+1}\rangle$.
Let
$\bar f_{ij} \in A_{\bs q,\Gr}$
be elements of Lemma \ref{lem proJ}.
Let
$h_{ij}$
be coefficients of the polynomials
$h_1,\dots,h_{N+1}$ in Theorem \ref{thm mV}.

\begin{thm}
[\cite{BMV}]
The map $\bar f_{ij} \mapsto h_{ij}$, \ $i=1,\dots, N+1$,
$j=1,\dots,d_i,\ d_i-j\not\in P$, extends uniquely to an algebra
isomorphism
$A_{\bs q,\Gr} \to A_{\T,\Phi}$.
\qed
\end{thm}

\begin{cor}
\label{cor bmv}
The elements
$h_{ij},\ i=1,\dots,{N+1},\ j=1,\dots,d_i,\ d_i-j\not\in P$,
generate $A_{\T,\Phi}$.
\qed
\end{cor}

\subsection{Proof of Theorem \ref{thm generate}}
\label{PRoofs}

By Lemma \ref{lem Generate} the subalgebra
of $A_{\T,\Phi}$ generated by all the elements $\bar G_{ij}$ contains all
the coefficients
$h_{ij}$. By Corollary \ref{cor bmv} the coefficients $h_{ij}$ generate
$A_{\T,\Phi}$. Theorem \ref{thm generate} is proved.

\section{Bethe algebra}
\label{sec algebra A_B}

\subsection{Lie algebra $\glnt$}

Let $\glnt=\gln\otimes\C[t]$ be the Lie algebra of $\gln$-valued
polynomials with the pointwise commutator. For $g\in\gln$, we set
$g(u)=\sum_{s=0}^\infty (g\otimes t^s)u^{-s-1}$.

We identify $\gln$ with the subalgebra $\gln\otimes1$ of constant polynomials
in $\glnt$. Hence, any $\glnt$-module has a canonical structure of
a $\gln$-module.

For each $a\in\C$, there exists an automorphism $\rho_a$ of $\glnt$,
\;$\rho_a:g(u)\mapsto g(u-a)$. Given a $\glnt$-module $M$, we denote by $M(a)$
the pull-back of $M$ through the automorphism $\rho_a$. As $\gln$-modules,
$M$ and $M(a)$ are isomorphic by the identity map.

We have the evaluation homomorphism,
${\glnt\to\gln}$, \;${g(u) \mapsto g u^{-1}}$.
Its restriction to the subalgebra $\gln\subset\glnt$ is the identity map.
For any $\gln$-module $M$, we denote by the same letter the $\glnt$-module,
obtained by pulling $M$ back through the evaluation homomorphism.

\subsection{Definition of row determinant}
\label{def rdet sec}
Given an algebra $A$ and an ${(N+1)\times (N+1)}$-matrix $C=(c_{ij})$ with entries in $A$,
we define its {\it row determinant\/} to be
\be
\rdet C\,=
\sum_{\;\si\in \Sigma_{N+1}\!} (-1)^\si\,c_{1\si(1)}c_{2\si(2)}\ldots c_{N+1\,\si(N+1)}\,.
\ee

\subsection{Definition of Bethe algebra}
\label{secbethe}

Define the { universal differential operator\/} $\D_\B$
by the formula
\be
\D_\B=\,\rdet\left( \begin{matrix}
\der_u-e_{11}(u) & - e_{21}(u)& \dots & -e_{N+1\,1}(u)\\[3pt]
-e_{12}(u) &\der_u-e_{22}(u)& \dots & -e_{N+1\,2}(u)\\[1pt]
\dots & \dots &\dots &\dots \\[1pt]
-e_{1\,N+1}(u) & -e_{2\,N+1}(u)& \dots & \der_u-e_{N+1\,N+1}(u)
\end{matrix}\right).
\ee
We have
\beq
\label{DK}
\D_\B=\,\der_u^{N+1}+\sum_{i=1}^{N+1}\,B_i\,\der_u^{{N+1}-i},
\qquad
B_i\,=\,\sum_{j=i}^\infty B_{ij} u^{-j}\,,
\qquad
B_{ij}\in \Uglnt\ .
\eeq
The unital subalgebra of $\Uglnt$ generated by $B_{ij}$, \,$i=1\lc {N+1}$,
\,$j\geq i$, is called the {\it Bethe algebra\/} and denoted by $\B$.

By \cite{T}, cf. \cite{MTV2}, the algebra $\B$ is commutative,
and $\B$ commutes with the subalgebra $U(\gln)\subset \Uglnt$.

As a subalgebra of $\Uglnt$, the algebra $\B$ acts on any $\glnt$-module
$M$. Since $\B$ commutes with $U(\gln)$, it preserves the $\gln$ weight subspaces
of $M$ and the subspace
$\sing M$ of $\gln$-singular vectors.

If $L$ is a $\B$-module, then the image of $\B$ in $\End(L)$ is called
the {\it Bethe algebra of\/} $L$.

\subsection{Bethe algebra of $\snash$}
\label{sec Bethe algebra}

For a partition $\bs\la$ with at most $N+1$ parts denote by $L_{\bs\la}$ the irreducible
$\gln$-module with highest weight $\bs\la$.

Let $\bs \La = (\bs\la^{(1)},\dots, \bs\la^{(n)})$, $\bs\la^{(\infty)}$,
$\bs z=(z_1,\dots,z_n)$ be partitions and numbers
as in Section \ref{subsubsection on master function}.
Denote
$
L_{\bs \La}
\ = \ L_{\bs\la^{(1)}}\otimes \dots \otimes L_{\bs\la^{(n)}}\,.
$
Let
\begin{gather*}
\nash\,=\,\{ v\in L_{\bs \La}\ | \
e_{ii}v= \la^{(\infty)}_iv\ {\rm for}\ i=1,\dots,N+1\} ,
\\
\snash\,=\, \{ v\in \nash\ | \
e_{ij}v= 0\ {\rm for}\ i<j \}
\end{gather*}
be the subspace of vectors of $\gln$-weight $\bs\la^{(\infty)}$ and
the subspace of $\gln$-singular vectors of $\gln$-weight $\bs\la^{(\infty)}$, respectively.
Consider on $\LL$ the $\glnt$-module structure of
the tensor product of evaluation modules,
$\LL = \otimes_{s=1}^n L_{\bs\la^{(s)}}(z_s)$.
Then $\snash$ is a $\B$-submodule.
We denote by $\BL$ the Bethe algebra of $\snash$.

\subsection{Shapovalov Form}
\label{Shapovalov Form}
Let $\tau : \gln \to \gln $ be the anti-involution sending
$e_{ij}$ to $e_{ji}$ for all $(i,j)$.
Let $M$ be a highest weight $\gln$-module with a
highest weight vector $w$.
{ The Shapovalov form} $S$ on $M$ is the unique
symmetric bilinear form such that
\be
S(w, w) = 1 ,
\qquad
S(xu, v) = S(u, \tau(x)v)
\ee
for all $u,v \in M$ and $x \in \gln$.

\medskip

Fix highest weight vectors $v_{\bs\la^{(s)}} \in L_{\bs\la^{(s)}}$,
$s=1,\dots,n$.
Define
a symmetric bilinear form on the tensor product
$
L_{\bs \La}
\ = \ L_{\bs\la^{(1)}}\otimes \dots \otimes L_{\bs\la^{(n)}}\,
$ by the formula
\beq\label{shap}
S_{\bs\La}\ =\ S_1 \otimes \cdots \otimes S_n ,
\eeq
where $S_s$ is the Shapovalov form on $L_{\bs\la^{(s)}}$.
The form $S_{\bs\la}$ is called the tensor Shapovalov form.

\medskip

\begin{thm}
[\cite{MTV2}]
\label{thm B has symm opers}
Consider $\LL$ as the $\glnt$-module $\otimes_{s=1}^n
L_{\bs\la^{(s)}}(z_s)$. Then any element $B\in \B$ acts on $\LL$ as a
symmetric operator with respect to the tensor Shapovalov form,
$S_{\bs\La}(Bu,v)= S_{\bs\La}(u,Bv)$ for any $u,v\in \LL$.
\qed
\end{thm}

\section{Weight function}
\label{sec weighT}

\subsection{Definition of the weight function}
\label{Defenition of a weight function}

Let $\bs \La = (\bs\la^{(1)},\dots, \bs\la^{(n)})$, $\bs\la^{(\infty)}$,
$\bs z=(z_1,\dots,$ $z_n)$ be partitions and numbers
as in Section \ref{subsubsection on master function}.
Recall the construction of a rational map
\be
\omega \ :\ \C^l \ \to \nash
\ee
called { the weight function}, see \cite{SV}, cf. \cite{M}, \cite{RSV}.

Denote by $P(\bs l,n)$ the set of sequences
$C\ = \ (c_1^1, \dots , c^1_{b_1};\ \dots ;\ c^n_1, \dots , c^n_{b_n})$ of integers from
$\{1, \dots , N\}$
such that for every
$i = 1, \dots , N$, the integer $i$ appears in $C$ precisely $l_i$ times.

Denote by $\Sigma(C)$ the set of all bijections $\sigma$
of the set
$\{1,\dots,l\}$ onto the set of variables
$\{t^{(1)}_1,\dots,t^{(1)}_{l_1},t^{(2)}_{1},\dots,t^{(2)}_{l_2},\dots,
t^{(N)}_{1},\dots,t^{(N)}_{l_{N}}\}$ with the following property.
For every $a=1,\dots,l$
the $a$-th element of the sequence $C$
equals $i$, if
$\sigma (a) = t^{(i)}_j$.

To every $C\in P(\bs l, n)$ we assign a vector
\be
e_Cv\ =\ e_{c_1^1+1,c_1^1} \dots e_{c_{b_1}^1+1,c_{b_1}^1}v_{\bs\la^{(1)}}
\otimes \cdots
\otimes e_{c_1^n+1,c_{1}^n} \dots e_{c_{b_n}^n+1,c_{b_n}^n}v_{\bs\la^{(n)}} \quad \in \ \nash .
\ee
To every $C\in P(\bs l, n)$ and
$\sigma\in \Sigma(C)$,
we assign a rational function
\be
\omega_{C,\sigma} \ =\ \omega_{\sigma;1,2,\ldots,b_1}(z_1)\
\cdots\
\omega_{\sigma; b_1+\dots + b_{n-1}+1,b_1+\dots + b_{n-1}+2,
\ldots,b_1+\dots + b_{n-1}+b_n}(z_n) ,
\ee
where
\be
\omega_{\sigma; a,a+1,\ldots, a+j}(z) \ =\ \frac 1 {(\sigma(a) - \sigma(a+1))
\ldots (\sigma(a+j-1) - \sigma(a+j))(\sigma(a+j) - z)} .
\ee
We set
\beq
\label{bethe vector}
\omega(\bs t) = \sum_{C\in P(\bs l,n)} \sum_{\sigma\in \Sigma (C)}
\omega_{C,\sigma} e_C v\ .
\eeq

\noindent
{\bf Examples.} If $n=2$ and $(l_1,l_2,\dots, l_N) = (1, 1, 0, \dots , 0)$, then
\begin{align*}
\omega(\bs t)\,={}&\,
\frac 1{(t^{(1)}_1-t^{(2)}_1)(t^{(2)}_1-z_1)}\;
e_{21}e_{32}v_{\bs\la^{(1)}}\otimes v_{\bs\la^{(2)}} +
\frac 1{(t^{(2)}_1-t^{(1)}_1)(t^{(1)}_1-z_1)}\;
e_{32}e_{21}v_{\bs\la^{(1)}}\otimes v_{\bs\la^{(2)}}
\\
{}+{}&\,
\frac 1{(t^{(1)}_1-z_1)(t^{(2)}_1-z_2)}\;
e_{21}v_{\bs\la^{(1)}}\otimes e_{32}v_{\bs\la^{(2)}} +
\frac 1{(t^{(2)}_1-z_1)(t^{(1)}_1-z_2)}\;
e_{32}v_{\bs\la^{(1)}}\otimes e_{21}v_{\bs\la^{(2)}}
\\
{}+{}&\,
\frac 1{(t^{(1)}_1-t^{(2)}_1)(t^{(2)}_1-z_2)}\;
v_{\bs\la^{(1)}}\otimes e_{21}e_{32}v_{\bs\la^{(2)}} +
\frac 1{(t^{(2)}_1-t^{(1)}_1)(t^{(1)}_1-z_2)}\;
v_{\bs\la^{(1)}}\otimes e_{32}e_{21}v_{\bs\la^{(2)}}\,.\!
\end{align*}
If $n=2$ and $(l_1,l_2,\dots, l_N) = (2, 0, \dots , 0)$, then
\begin{align*}
\omega(\bs t)\,={}& \,
(\frac 1{(t_1^{(1)}-t_2^{(1)})(t^{(1)}_2-z_1)}+
\frac 1{(t^{(1)}_2-t^{(1)}_1)(t^{(1)}_1-z_1)})\;
e_{21}^2v_{\bs\la^{(1)}}\otimes v_{\bs\la^{(2)}}
\\
{}+{}&\,
(\frac 1{(t^{(1)}_1-z_1)(t^{(1)}_2-z_2)}+
\frac 1{(t^{(1)}_2-z_1)(t^{(1)}_1-z_2)})\;
e_{21}v_{\bs\la^{(1)}}\otimes e_{21}v_{\bs\la^{(2)}}
\\
{}+{}&\,
(\frac 1{(t^{(1)}_1-t^{(1)}_2)(t^{(1)}_2-z_2)}+
\frac 1{(t^{(1)}_2-t^{(1)}_1)(t^{(1)}_1-z_2)})\;
v_{\bs\la^{(1)}}\otimes e_{21}^2v_{\bs\la^{(2)}}\,.
\end{align*}

\begin{lem} [Lemma 2.1 in \cite {MV2}]
\label{well def}
The weight function is regular on $U$.
\qed
\end{lem}

\subsection{Grothendieck residue and Hessian}
\label{Grothendieck residue}

Let
\be
{\rm Hess}\, \log\,\Phi \ = \ \det \left(
\frac{\partial^2}{\partial t^{(i)}_j
\partial t^{(i')}_{j'}} \log\, \Phi
\right)
\ee
be the Hessian of $\log\,\Phi$.
Let $\T\in U$ be a critical point of $\Phi$.
Denote by $H_\T$ the image of the Hessian in
the local algebra
$\AT$. It is known that $H_\T$ is nonzero and
the one-dimensional subspace $\C H_\T\subset \AT$ is the annihilating ideal of the maximal ideal
$\frak m_\T \subset \AT$.

Let
$\rho_{\T} : A_{{\T},\Phi} \to \C$, be the
Grothendieck residue,
\be
f \ \mapsto \ \frac 1{(2\pi i)^l}\,\Res_{\T}
\frac{ f}{\prod_{ij} \Psi_{ij}}\ .
\ee
It is known that
$
\rho_{\T}(H_\T) = \mu_\T,
$
where $\mu_\T = \dim\, A_{\T,\Phi}$ is the Milnor number of the critical point
$\T$.
Let $(\,,\,)_{\T}$ be the Grothendieck residue
bilinear form on $A_{\T, \Phi}$\,,
\be
( f, g)_{\T}\ =\ \rho_{\T} (f g)\ .
\ee
It is known that $(\,,\,)_\T$ is nondegenerate. These facts see for example in
Section~5 of~\cite{AGV}.

\subsection{Projection of the weight function}

Let $\T\in C_\Phi$ be a critical point of $\Phi$. Let
$$
\omega_\T \in \nash \otimes A_{\T, \Phi}
$$
be the element induced by the weight function.
Let $S_{\bs \La}$ be the tensor Shapovalov form on
$L_{\bs\La}$.

\begin{thm} [\cite{MV2}, \cite{V}]
\label{thm MV2, V}
We have
\beq
\label{norm formula}
S_{\bs \La}(\omega_\T,\omega_\T) =
H_\T\ .
\eeq
\qed
\end{thm}

\begin{thm}
[\cite{SV}]
\label{thm [SV]}
The element $\omega_\T$ belongs to
$\snash\otimes \AT $
where
\\
$\snash \subset \nash$ is the subspace of singular vectors.
\qed
\end{thm}

Theorem \ref{thm [SV]} is a direct corollary of Theorem 6.16.2 in \cite{SV}, see also
\cite{RV} and \cite{B}.

\subsection{Bethe ansatz}
Let $\T\in C_\Phi$ be a critical point of $\Phi$.
Consider the differential operator
\be
\label{formula htu}
\D_{\Phi}
\ =\
\der_u^{N+1}+\sum_{i=1}^{N+1}\,G_i\,\der_u^{N+1-i},
\qquad
G_i = \sum_{j=i}^\infty G_{ij} u^{-j} ,
\qquad
G_{ij} \in \C[\bs t] ,
\ee
described by \Ref{formula h(t,u)}, and
projections $\bar G_{ij}$ of its coefficients to $\AT$.
Consider the differential operator
\be
\label{diff op}
\D_\B=\,\der_u^{N+1}+\sum_{i=1}^{N+1}\,B_i\,\der_u^{{N+1}-i},
\qquad
B_i\,=\,\sum_{j=i}^\infty B_{ij} u^{-j}\,,
\qquad
B_{ij}\in \Uglnt\ ,
\ee
described by \Ref{DK}.

\begin{thm} [\cite {MTV2}]
\label{thm MtV}
For any $i=1,\dots,N+1,\,j\geq i$, we have
\beq
\label{bethe formula}
(B_{ij}\otimes 1)\, \omega_\T
\ = \
(1\otimes \bar G_{ij})\, \omega_\T
\eeq
in $\snash\otimes \AT$.
\qed
\end{thm}

This statement is the Bethe ansatz method to construct eigenvectors of
the Bethe algebra in the $\gln$ Gaudin model starting with a critical
point of the master function.

\subsection{Main result}
\label{sec Main result}

Let $g_1,\dots, g_{\mu_\T}$ be a basis of $\AT$
considered as a $\C$-vector space. Write
$\omega_\T =\sum_i v_i
\otimes g_i$, with $ v_i \in \snash$.
Denote by $M_\T \subset {\rm Sing}\,\nash$ the vector subspace spanned by
$v_1,\dots,v_{\mu_\T}$.
Define a linear map
\beq
\label{map psi}
\al\ :\ A_{\T, \Phi}\ \to\ M_\T\ ,
\qquad
f\ \mapsto\ (f,\omega_\T)_\T = \sum_{i=1}^{\mu_\T}\ (f,g_i)_\T\,v_i\ .
\eeq

\begin{thm}
\label{thm main}
Let $\T\in C_\Phi$. Then the following statements hold:
\begin{enumerate}
\item[(i)]
The subspace $M_\T\subset \snash$ is a $\B$-submodule.
Let $\AB \subset \End\,(M_\T)$ be the Bethe algebra of $M_\T$.
Denote by $\bar B_{ij}$ the image in $\AB$ of generators $B_{ij}\in \B$.

\item[(ii)]
The map $\al : A_{\T, \Phi} \to M_\T$ is an isomorphism of vector spaces.
\item[(iii)]
The map $\bar G_{ij} \mapsto \bar B_{ij}$ extends uniquely to an algebra
isomorphism $\beta: \AT \to \AB$.
\item[(iv)]

The isomorphisms $\al$ and $\beta$ identify the regular representation
of $\AT$ and the $\B$-module $M_\T$\,,\ that is, for any $f,g\in\AT$ we have
$\al(fg)=\beta(f)\al(g)$.
\end{enumerate}

\end{thm}

\begin{cor}
\label{cor value at crit point}
Let $\T\in C_\Phi$.
Then the value $\omega(\T)$ of the weight function at $\T$ is a nonzero vector of
$\snash$.
\end{cor}

Theorem \ref{thm main} and Corollary \ref{cor value at crit point}
are proved in Section \ref{prooFs}.

\section{Proof of Theorem \ref{thm main}}
\label{pRoofs}

\subsection{Proof of part (i) of Theorem \ref{thm main}}
\label{lem M B-submodule}
It is enough to show that for any $f\in \AT$ and any $(i,j)$ we have
$B_{ij}\al(f)\in M_\T$. Indeed, we have
\beq
\label{useful}
B_{ij}\al(f) = \sum_{l=1}^{\mu_\T} (f,g_l)_\T B_{ij}v_l =
\sum_{l=1}^{\mu_\T} (f,\bar G_{ij}g_l)_\T\bar v_l=
\sum_{l=1}^{\mu_\T} (\bar G_{ij}f,g_l)_\T\bar v_l =
\al(\bar G_{ij}f) .
\eeq
Here the second equality follows from Theorem \ref{thm MtV}
and the third equality follows from properties of the
Grothendieck residue form.

\subsection{Bilinear form $(\,,\,)_S$}

Define a symmetric bilinear form $(\,,\,)_S$ on $\AT$,
\be
(f,g)_S = S_{\bs \La}(\al (f),\al (g)) =
\sum_{i,j=1}^{\mu_\T}\, S_{\bs \La}(v_i,v_j)\, (f,g_i)_\T\,(g,g_j)_\T\
\ee
for all $f,g \in \AT$.

\begin{lem}
\label{lem property ()S}
For all $ f,g, h \in \AT$ we have
$(fg, h)_S \ =\ (f, g h)_S$.
\end{lem}

\begin{proof}
By Theorem \ref{thm generate} the elements $\bar G_{ij}$ generate $\AT$.
We have
$(\bar G_{ij}f,h)_S =$
\\
$ S_{\bs \La}(\al(\bar G_{ij}f),\al(h)) =
S_{\bs \La}(\bar B_{ij}\al (f),\al (h)) =
S_{\bs \La}(\al (f),\bar B_{ij}\al (h)) =
S_{\bs \La}(\al (f),\al (\bar G_{ij}h)) =
(f,\bar G_{ij}h)_S.$
Here the third equality follows from Theorem \ref{thm B has symm opers}.
\end{proof}

\begin{lem}
\label{lem F}
There exists $F\in \AT$ such that
$(f,h)_S = (Ff,h)_\T$
for all $ f, h \in \AT$.
\end{lem}

\begin{proof}
Consider the linear function $\AT \to \C,\,
h \mapsto (1,h)_S$. The form $(\,,\,)_\T$ is nondegenerate. Hence there exits
$F\in \AT$ such that $(1,h)_S = (F,h)_\T$ for all $h\in \AT$.
Now the lemma follows from Lemma \ref{lem property ()S}.
\end{proof}

\subsection{Auxiliary lemmas}

\begin{lem}
\label{lemma 4}

For any $f\in \AT$, we have
\beq
\label{123}
f H_\T = \frac {1}{\mu_\T} (f, H_\T)_\T H_\T\, .
\eeq

\end{lem}

\begin{proof}
The lemma follows from the fact that formula \Ref{123} evidently
holds for $1 \in \AT$ and for
any element of the maximal ideal.
\end{proof}

For $f\in\AT$, denote by $L_f$ the linear operator $\AT\to\AT,\, h\mapsto fh$.

\begin{lem}
\label{lemma 5}
We have\
${\rm tr}\,L_f\, =\,(f, H_\T)_\T .$
\end{lem}

\begin{proof}
The linear function $\AT \to \C,\ f\mapsto {\rm tr}\,L_f,$ is
such that $1 \mapsto \mu_\T$ and $f\mapsto 0$ for all $f\in \frak m_\T$.
Hence this function equals the linear function $f \mapsto (f,\,H_\T)_\T$.
\end{proof}

Let $g_1^*,\dots,g_{\mu_\T}^*$ be the basis of $\AT$ dual to the basis
$g_1,\dots,g_{\mu_\T}$ with respect to the form $(\,,\,)_\T$. Then
$H_\T=\sum_{i=1}^{\mu_\T} (H_\T,g_i^*)_\T g_i$.
Indeed for any $f\in \AT$, we have
$f = \sum_i\, (f,g_i^*)_\T g_i$.

\begin{lem}
\label{new lem}
We have $\sum_{i=1}^{\mu_\T} g_i^*g_i = H_\T$ .
\end{lem}

\begin{proof}
For $f\in \AT$, we have\ {}
${\rm tr}\, L_f\, =\, \sum_i(g_i^*, f g_i)_\T =
(\sum_i\,g_i^*g_i, f)_\T .
$
By Lemma \ref{lemma 5}, we get \ {}
$
(\sum_i\,g_i^*g_i, f)_\T\, =\, (H_\T, f)_\T .
$\
Hence $\sum_i\,g_i^*g_i = H_\T$, since the form
$(\,,\,)_\T$ is nondegenerate.
\end{proof}

\begin{lem}

Let $F\in \AT$ be the element defined in Lemma \ref{lem F}.
Then $F$ is invertible, $F H_\T = H_\T$, and
the form $(\,,\,)_S$ is nondegenerate.
\end{lem}

\begin{proof} By definitions we have
\be
(f,h)_S = \sum_{ij} S_{\bs \La}(v_i,v_j) (g_i,f)_\T (g_j,h)_\T\
\ee
and
\be
(f,h)_S = (Ff,h)_\T = \sum_{i} (g_i,Ff)_\T (g_i^*,h)_\T =
\sum_{i} (F g_i,\,f)_\T (g_i^*,h)_\T\, .
\ee
Hence
$
\sum_{ij} S_{\bs \La}(v_i,v_j) g_i \otimes g_j =
\sum_{i} F g_i \otimes g_i^*
$
and therefore by Lemma \ref{new lem} we get
\be
\sum_{ij} S_{\bs \La}(v_i,v_j) g_i g_j = \sum_{i} F g_i g_i^* = FH_\T .
\ee
By Theorem \ref{thm MV2, V}, \
$
\sum_{ij} S_{\bs \La}(v_i,v_j)g_i g_j = H_\T .
$
Hence $FH_\T = H_\T$, the element
$F$ is invertible, and
the form $(\,,\,)_S$ is nondegenerate.
\end{proof}

\subsection{Proof of Theorem \ref{thm main} and
Corollary \ref{cor value at crit point}}
\label{prooFs}
Part (i) of Theorem \ref{thm main}
is proved in Section \ref{lem M B-submodule}.

Assume that $\sum_{i=1}^{\mu_\T} \la_iv_i = 0$. Denote $h = \sum_i \la_i
g_i^*$. Then $\al(h) = 0$ and $(f,h)_S = S_{\bs \La}(\al(f),
\al(h)) = 0$ for all $f\in\AT$. Hence $h=0$ since $(\,,\,)_S$ is
nondegenerate. Therefore, $\la_i = 0$ for all $i$ and the vectors
$v_1,\dots,v_{\mu_\T}$ are linearly independent. We have $\al(g_i^*) = v_i$
for all $i$.
That proves part (ii) of Theorem \ref{thm main}.

Parts (iii-iv) easily follow from part (ii) and formula \Ref{useful}.

We have
\beq
\label{eigenvectoR}
\mu_\T \,\omega(\T) = (H_\T, \omega_\T)_\T\ = \al(H_\T) .
\eeq
That implies that $\omega(\bs\T)$ is a nonzero vector.

\section{Concluding remarks}
\label{Concluding remarks}
\begin{thm}
\label{thm many point}
Let $C_\Phi = \{\T_1,\dots,\T_k\}$, be the critical set of $\Phi$ in $U$.
Let $M_{\T_s}\subset \snash$, $s=1,\dots,k$, be the corresponding subspaces
defined in Section \ref{sec Main result}. Then the sum of these subspaces is
direct.
\end{thm}

\begin{proof}
It follows from Theorem \ref{thm main} that for any $s$ and any $(i,j)$ the
operator $\bar B_{ij} - G_{ij}(\T_s)$ restricted to $M_{\T_s}$ is nilpotent.
Moreover, the differential operators $\D_\Phi|_{\bs t =\T_s}$, \,$s=1,\dots,k$,
which contain eigenvalues of the operators $B_{ij}$, are distinct.
These observations imply Theorem~\ref{thm many point}.
\end{proof}

%
%
%

Let \,$\al(A_\Phi) = \oplus_{s=1}^k M_{p_s}$. Denote
by $A_B$ the image of $\B$ in $\End(\al(A_\Phi))$.
Consider the isomorphisms
\be
\al =\oplus_{s=1}^k\al_s\ :\ \oplus_{s=1}^kA_{\T_s, \Phi}\ \to\ \oplus_{s=1}^kM_{\T_s} \,,
\qquad
\beta=\oplus_{s=1}^k\beta_s \ :\ \oplus_{s=1}^kA_{\T_s, \Phi}\ \to\ \oplus_{s=1}^k A_{\T_s,B}\,
\ee
of Theorem \ref{thm main}.

\begin{cor}
\label{cor main}
We have
\begin{enumerate}
\item[(i)]
$A_B=\oplus_{s=1}^k A_{p_s,B}$;
\item[(ii)] The isomorphisms
$\al$, $\beta$
identify the regular representation
of the algebra $A_\Phi$ and the $A_B$-module $\al(A_\Phi)$.
\end{enumerate}
\end{cor}

The corollary follows from Theorems \ref{thm main} and \ref{thm many point}.

\end{document}